\documentclass[12pt]{article}
\usepackage{amsmath,amsfonts,amssymb,amsthm}
\newcommand\bi[2]{{{#1}\atopwithdelims(){#2}}}
\newcommand\qbi[3]{{{#1}\atopwithdelims[]{#2}}_{#3}}

\usepackage{epsfig} \usepackage{amssymb,amsfonts,latexsym}
\newtheorem{theo}{Theorem}
 \newtheorem{coro}{Corollary}
 \newtheorem{lem}{Lemma}

\newenvironment{rem}{\noindent{\bf Remark. }}{\smallskip}
\newcommand{\la}{\lambda}
\title{A generalization  of Kawanaka's identity for
 Hall-Littlewood polynomials and applications }
\author{Masao Ishikawa$^1$, Fr\'ed\'eric Jouhet$^2$ and Jiang Zeng$^{2}$}
\date{}
\begin{document}
\maketitle
\begin{center}
 \small $^1$ Department of Mathematics, Faculty of Education\\
Tottori University, Tottori 680-8551, Japan,\\
\small \texttt{ishikawa@fed.tottori-u.ac.jp}\\
\vspace{10 pt}
 \small $^2$ Institut Camille Jordan,
Universit\'e Claude Bernard (Lyon I)\\
43, bd du 11 Novembre 1918, 69622 Villeurbanne Cedex, France \\
\small \texttt{\{jouhet,zeng\}@math.univ-lyon1.fr}\\
\end{center}
\begin{abstract}
An infinite summation formula of Hall-Littlewood polynomials due
to Kawanaka is generalized to a finite summation formula, which
implies, in particular, twelve more multiple $q$-identities of
Rogers-Ramanujan type than those previously
found by Stembridge and the last two authors.\\

\end{abstract}

\section{Introduction}
Recently, starting from two infinite summation formulae for
Hall-Littlewood polynomials,  two of the present authors~\cite{JZ2} have
generalized a method due to Macdonald~\cite{Ma}  to obtain new
finite summation formulae for these polynomials. This approach permits them
to extend Stembridge's list of multiple
$q$-series identities of Rogers-Ramanujan type~\cite{St}. Conversely these
symmetric functions identities can be viewed as a generalization
of Rogers-Ramanujan identities. In view of  the numerous formulae
of Rogers-Ramanujan type~\cite{Sl2} one may speculate that
there should be more such  generalizations starting from other infinite summation formulae for
Hall-Littlewood polynomials.
However, as pointed out in
\cite{JZ2}, when one passes from an infinite summation to a finite summation,
one may need to modify the
coefficients normalizing Hall-Littlewood polynomials
 in order to obtain some  useful formulae.

In this paper we take up on Kawanaka's  new
infinite summation identities of Hall-Littlewood polynomials~\cite{Kaw} and
 show that one of his formulae has a finite summation generalization.

We first need to recall some
standard notations of $q$-series, which can be found in~\cite{GR}.
Set $(x)_0:=(x;q)_0=1$ and for $n\geq 1$
$$
(x)_n:=(x;q)_n=\prod_{k=0}^n(1-xq^{k}),\quad  (x)_\infty:=(x;q)_\infty=\prod_{k=0}^\infty(1-xq^{k}).
$$
For $n\geq 0$ and $r\geq 1$, set
$$
(a_1,\cdots,a_r;\,q)_n=\prod_{i=1}^r(a_i)_n\,\,\qquad (a_1,\cdots
,a_r;\,q)_\infty=\prod_{i=1}^r(a_i)_\infty.
$$
The $q$-binomial identity~\cite{An76} then reads  as follows~:
\begin{equation}\label{qbinomial}
\sum_{m\geq
0}\frac{(a)_m}{(q)_m}x^m=\frac{(ax)_\infty}{(x)_\infty},
\end{equation}
which reduces to the finite $q$-binomial identity by substitution
$a\rightarrow q^{-n}$ and $x\rightarrow q^nx$:
\begin{equation}\label{qbin}
(x)_n=\sum_{m\geq 0}(-1)^mq^{\left({m\atop
2}\right)}\qbi{n}{m}{}x^m
\end{equation}
and to the following identity of Euler when  $a=0$~:
\begin{equation}\label{euler}
\frac{1}{(x)_\infty}=\sum_{m\geq 0}\frac{x^m}{(q)_m}.
\end{equation}

Let $n\geq 1$ be a fixed integer and  $S_n$ denote  the group of
permutations of the set $\{1,\, 2, \ldots, n\}$. Let
$X=\{x_1,\ldots, x_n\}$ be a set of indeterminates and $q$ a
parameter. For each \emph{partition} $\lambda=(\la_1,\ldots,
\la_n)$ of length $\leq n$, if $m_i:=m_i(\la)$ is  the
multiplicity of part $i$  in $\la$, then we also denote $\la$ by
$(1^{m_1}\,2^{m_2}\, \ldots)$. Recall that the Hall-Littlewood
polynomials $P_{\la}(X,q)$ are defined by~\cite[p.208]{Ma}~:
$$
P_{\la}(X,q)= \prod_{i\geq
1}\frac{(1-q)^{m_i}}{(q)_{m_i}}\,\sum_{w\in
S_n}w\left(x_1^{\la_1}\ldots x_n^{\la_n}
\prod_{i<j}\frac{x_i-qx_j}{ x_i-x_j}\right).
$$
Since (\cite[p. 207]{Ma})
$$
\sum_{w\in
S_n}w\left(\prod_{i<j}\frac{x_i-qx_j}{ x_i-x_j}\right)=\frac{(q)_n}{(1-q)^n},
$$
we see that the coefficient of $x_1^{\la_1}\ldots x_n^{\la_n}$ in
$P_\la$ is 1.
Set
$$
\Phi(X):=\prod_{i}\frac{1+qx_i}{1-x_i}\prod_{j<k}\frac{1-q^2x_jx_k}{1-x_jx_k}.
$$
Our starting point is  the following result due to Kawanaka~\cite{Kaw}:
\begin{equation}\label{kawa1}
 \sum_\la \left(\prod_{i\geq 1}(-q)_{m_i}\right)
 P_\la(X,q^2)=\Phi(X).
\end{equation}
Since Kawanaka's proof uses the representation theory of groups we shall give
another proof using Pieri's rule for Hall-Littlewood polynomials.

For each sequence $\xi\in\{\pm 1\}^n$, set $X^{\xi}:=\{x_1^{\xi_1},\ldots,x_n^{\xi_n}\}$.
Our  finite extension of Kawanaka's formula then reads as follows:
\begin{theo} For $k\geq 1$ the following identity holds
\begin{equation}\label{extkawa1}
\sum_{\la_1\leq
k}\left(\prod_{i=1}^{k-1}(-q)_{m_i}\right)P_\la(X,q^2)=\sum_{\xi\in \{\pm
1\}^n}\Phi(X^{\xi})\prod_ix_i^{k(1-\xi_i)/2}.
\end{equation}
\end{theo}
\begin{rem}  In the  case $q=0$,  the right-hand side of
(\ref{extkawa1}) can be written as a quotient of determinants and
the formula reduces to a known identity of Schur
functions~\cite{JZ1}.
\end{rem}

For any partition $\la$ it will be convenient  to adopt  the
following notation~: $$(x)_\la:= (x;\,
q)_\la=(x)_{\la_1-\la_2}(x)_{\la_2-\la_3}\cdots
$$
Note that  this is not the standard notation
for $(x)_\la$ and corresponds to $b_{\la'}(q)$ in~\cite[p.210]{Ma}.

We also introduce the following generalization of  $q$-binomial coefficients
$$
\qbi{n}{\la}{}:=\frac{(q)_n}{(q)_{n-\la_1}(q)_{\la}},
$$
with the
convention that $\qbi{n}{\la}{}=0$ if $\la_1>n$. If $\la=(\la_1)$
we recover the classical $q$-binomial coefficient. Finally, for
any partition $\la$ we denote by $l(\la)$ the length of $\la$,
i.e., the number of its positive parts, and define
$n(\la):=\sum_i\left({\la_i\atop 2}\right)$ and  $n_2(\la)=\sum_i \la_i^2$. When  $x_i=zq^{2i-2}$
for $i\geq 1$ and then $z$ is replaced by $zq$, formula (\ref{extkawa1}) specializes to the
following identity.
\begin{coro} For $k\geq 1$ there holds
\begin{eqnarray}
&&\hskip -1 cm\sum_{l(\la)\leq k}\left(\prod_{i=1}^{k-1}(-q)_{\la_i-\la_{i+1}}\right)z^{|\la|}q^{n_2(\la)}\left[{n\atop \la}\right]_{q^2}\label{csqextkawa1}\\
&&=\sum_{r=0}^n(-1)^rz^{kr}q^{(k+1)r^2}\left[{n\atop r}\right]_{q^2}
\frac{(-z)_{2n+1}}{(z^2q^{2r};q^2)_{n+1}}(1-zq^{2r}).\nonumber
\end{eqnarray}
\end{coro}

Now, as in \cite{JZ2, St}, we can prove the following key
$q$-identity which allows to produce identities of
Rogers-Ramanujan type :
\begin{theo} For $k\geq 1$,
\begin{eqnarray}
&&\sum_{l(\la)\leq k}z^{|\la|}q^{2n(\la)}\frac{(a,b;q^{-1})_{\la_1}}{(-q)_{\la_k}(q)_\la}\label{ab}\\
&&\hskip -0.5 cm=\frac{(-z/q)_{\infty}}{(abz)_\infty}\sum_{r\geq 0}(-1)^rz^{kr}q^{r+(2k+2)\bi{r}{2}}\frac{(a,b;q^{-1})_r}{(q^2;q^2)_r}\frac{(azq^r,bzq^r)_\infty}{(z^2q^{2r-2};q^2)_{\infty}}(1-zq^{2r-1}).\nonumber
\end{eqnarray}
\end{theo}

This paper is organized as follows. In Section~2 we give a new
proof of Kawanaka's formula using Pieri's rule for Hall-Littlewood
polynomials since Kawanaka's original proof uses the
representation theory of groups. In section~3, we derive from
Theorem~2 twelve
 multiple analogs of Rogers-Ramanujan type
identities. In section~4 we give the proofs of Theorem~1 and
Corollary~1, and some consequences, and defer the elementary
proof, i.e., without using the Hall-Littlewood polynomials, of
Theorem~2, Corollary~1 and other multiple $q$-series identities to
section~5. To prove Theorems 1 and 2  we apply the generating
function technique which was developped in \cite{JZ2, Ma, St}.

\section{Another proof of Kawanaka's identity}
Recall \cite[p.230, Ex.1]{Ma} the following summation of Hall-Littlewood polynomials :
$$
\sum_{\mu}P_{\mu}(X,q)
=\prod_{i}\frac1{1-x_i}\prod_{i<j}\frac{1-qx_ix_j}{1-x_ix_j}.
$$
By replacing $q$ by $q^2$, we get
\begin{equation}
\label{HL_sum} \sum_{\mu}P_{\mu}(X,q^2)
=\prod_{i}\frac1{1-x_i}\prod_{i<j}\frac{1-q^2x_ix_j}{1-x_ix_j}.
\end{equation}
Note that
\begin{equation}
\label{elementary} \sum_{r\geq0}e_{k}(X)q^k=\prod_{i}(1+qx_i),
\end{equation}
where $e_r(X)$ stands for the $r$-th elementary symmetric
function. Identities \thetag{\ref{HL_sum}} and \thetag{\ref{elementary}}
imply
\begin{equation*}
\sum_{\mu}\sum_{r}q^rP_{\mu}(X,q^2)e_r(X)
=\prod_{i}\frac{1+qx_i}{1-x_i}
\prod_{i<j}\frac{1-q^2x_ix_j}{1-x_ix_j}.
\end{equation*}
From \cite[p.209, (2.8)]{Ma}, we have
\begin{equation*}
P_{(1^r)}(X,q)=e_r(X),
\end{equation*}
and this shows that
\begin{equation*}
\sum_{\mu}\sum_{r}q^rP_{\mu}(X,q^2)P_{(1^r)}(X,q^2)
=\prod_{i}\frac{1+qx_i}{1-x_i}
\prod_{i<j}\frac{1-q^2x_ix_j}{1-x_ix_j}.
\end{equation*}
Let $f^{\lambda}_{\mu\nu}(q)$ be the coefficients defined by
\begin{equation*}
P_{\mu}(X,q)P_{\nu}(X,q)=\sum_{\lambda}f^{\lambda}_{\mu\nu}(q)
P_{\lambda}(X,q),
\end{equation*}
then, by \cite[p.215 (3.2)]{Ma} we have
\begin{equation*}
f^{\lambda}_{\mu(1^m)}(q)=\prod_{i\geq 1}
\left[{{\lambda_i'-\lambda_{i+1}'}\atop{\lambda_i'-\mu_i'}}\right]_{q}
\end{equation*}
(and therefore $f^{\lambda}_{\mu(1^m)}(q)=0$ unless $\lambda\setminus\mu$
is a vertical $m$-strip, or $m$-vs, which means $\lambda\subset\mu$, $|\lambda\setminus\mu|=m$ and there is at most one cell in each row of the Ferrers diagram of $\lambda\setminus\mu$). Thus we have
\begin{align*}
&\sum_{\lambda}\sum_{{\mu}\atop{\lambda\setminus\mu\;\text{vs}}} q^{|\lambda-\mu|}\prod_{i\geq1}
\left[{{\lambda_i'-\lambda_{i+1}'}\atop{\lambda_i'-\mu_i'}}\right]_{q^2}P_{\lambda}(X,q^2)\\
&\qquad\qquad=\prod_{i}\frac{1+qx_i}{1-x_i}
\prod_{i<j}\frac{1-q^2x_ix_j}{1-x_ix_j}.
\end{align*}
Applying the identity~(see \cite{An76}, and \cite{Ze} for a bijective
proof) :
\begin{equation}\label{csqeuler}
\sum_{k=0}^{n}q^{k}\left[{{n}\atop{k}}\right]_{q^2}
=\prod_{k=1}^{n}(1+q^{k}),
\end{equation}
we conclude that
$$
\sum_{{\mu}\atop{\lambda\setminus\mu\;\text{vs}}}
q^{|\lambda-\mu|}\prod_{i\geq 1}
\left[{{\lambda_i'-\lambda_{i+1}'}
\atop{\lambda_i'-\mu_i'}}\right]_{q^2} =\prod_{i\geq
1}\prod_{k=1}^{\lambda_i'-\lambda_{i+1}'}(1+q^k),
$$
which is precisely what we desired to prove.

\begin{rem}
For a node $v=(i,j)$ in the diagram of $\la$, the arm-length $a(v)$ and the leg-length $l(v)$ of $\la$ at $v$ are defined by $a(v)=\la_i-j$ and $l(v)=\la'_j-i$ respectively. Kawanaka~\cite[(5.2)]{Kaw} proved
another identity for Hall-Littlewood polynomials~:
\begin{equation}
\sum_{\la}q^{o(\la)/2}\left(\prod_{{v\in\lambda,\;a(v)=0 \atop l(v) \text{even}}}(1-q^{l(v)+1})\right)P_\la(X,q)=\prod_{i\leq
j}\frac{1-qx_ix_j}{1-x_ix_j},
\end{equation}
where the sum on the left is taken over all partitions $\la$ such
that $m_i(\la)$ is even for odd $i$ and
$$
o(\la)=\sum_{i\;\textrm{odd}}m_i(\la).
$$
It would be possible to prove this identity in the same manner as
above.

 There is a related
identity about Hall-Littlewood polynomials in Macdonald's
book~\cite[p. 219]{Ma}~:
\begin{equation}
\sum_\la q^{n(\la)}\left(\prod_{j=1}^{l(\la)}(1+q^{1-j}y)\right)P_\la(X,q) =\prod_{i\geq 1}\frac{1+x_iy}{1-x_i}.
\end{equation}
\end{rem}
\section{Multiple identities of Rogers-Ramanujan type}
We shall derive several identities  of Rogers-Ramanujan type from
Theorem~2.
First we note that if $z=q^2$ identity (\ref{ab}) reduces to
\begin{eqnarray}
&&\sum_{l(\la)\leq k}q^{|\la|+n_2(\la)}\frac{(a,b;q^{-1})_{\la_1}}{(-q)_{\la_k}(q)_\la}
={1\over (q,abq^2)_\infty}\label{zq2}\\
&&\times\sum_{r\geq
0}(-1)^rq^{(2k+1)r+(2k+2)\bi{r}{2}}(a,b;q^{-1})_r(aq^{r+2},bq^{r+2})_\infty
(1-q^{2r+1}),\nonumber
\end{eqnarray}
and if $z=q$ it becomes
\begin{eqnarray}
&&\sum_{l(\la)\leq k}q^{n_2(\la)}\frac{(a,b;q^{-1})_{\la_1}}{(-q)_{\la_k}(q)_\la}
={1\over (q,abq)_\infty} \label{zq}\\
&&\times\left((aq,bq)_\infty+2\sum_{r\geq
1}(-1)^rq^{(k+1)r^2}(a,b;q^{-1})_r(aq^{r+1},bq^{r+1})_\infty\right).\nonumber
\end{eqnarray}

We need the following two forms of {\em Jacobi triple product}
identity~\cite[p.21]{An76}~:
\begin{eqnarray}
J(x,\, q):=(q,\, x,\,  q/x)_\infty
&=&\sum_{r= 0}^\infty (-1)^rx^rq^{\left({r \atop 2}\right)}(1-q^{2r+1}/x^{2r+1})\label{Ja1}\\
&=&1+\sum_{r= 1}^\infty (-1)^rx^rq^{\left({r \atop
2}\right)}(1+q^r/x^{2r})\label{Ja2}.
\end{eqnarray}

\begin{theo} For $k\geq 1$, the following identities hold
\begin{eqnarray}
\sum_{l(\la)\leq k}\frac{q^{|\la|+n_2(\la)}}{(-q)_{\la_k}(q)_\la}
&=&\frac{(q^{2k+2},q^{2k+1},q;q^{2k+2})_\infty}{(q)_\infty},\label{krr1}\\
\sum_{l(\la)\leq
k}\frac{q^{|\la|+n_2(\la)-(\la_1^2+\la_1)/2}(-q)_{\la_1}}{(-q)_{\la_k}(q)_\la}
&=&\frac{(-q)_\infty}{(q)_\infty}(q^{2k+1},q^{2k},q;q^{2k+1})_\infty\label{krr2},\\
\sum_{l(\la)\leq k}
\frac{q^{2|\la|+2n_2(\la)-\la_1^2}(-q;q^2)_{\la_1}}{(-q^2;q^2)_{\la_k}(q^2;q^2)_\la}
&=&\frac{(-q;q^2)_\infty}{(q^2;q^2)_\infty}(q^{4k+2},q^{4k+1},q;q^{4k+2})_\infty\label{krr3},\\
\sum_{l(\la)\leq k}
\frac{q^{2|\la|+2n_2(\la)-2\la_1^2-\la_1}(-q)_{2\la_1}}
{(-q^2;q^2)_{\la_k}(q^2;q^2)_\la}
&=&\frac{(-q)_\infty}{(q)_\infty}(q^{4k},q^{4k-1},q;q^{4k})_\infty,\label{krr4}
\end{eqnarray}
\begin{eqnarray}
\sum_{l(\la)\leq k}\frac{q^{n_2(\la)}}{(-q)_{\la_k}(q)_\la}
&=&\frac{(q^{2k+2},q^{k+1},q^{k+1};q^{2k+2})_\infty}{(q)_\infty}\label{krr9},\\
\sum_{l(\la)\leq
k}\frac{q^{n_2(\la)-(\la_1^2+\la_1)/2}(-q)_{\la_1}}{(-q)_{\la_k}(q)_\la}
&=&\frac{(-1)_\infty}{(q)_\infty}(q^{2k+1},q^{k},q^{k+1};q^{2k+1})_\infty\label{krr11},\\
\sum_{l(\la)\leq k}\frac{q^{2n_2(\la)-\la_1^2}
(-q;q^2)_{\la_1}}{(-q^2;q^2)_{\la_k}(q^2;q^2)_\la}
&=&\nonumber\\
&&\hskip -1 cm\frac{(-q;q^2)_\infty}{(q^2;q^2)_\infty}(q^{4k+2},q^{2k+1},q^{2k+1};q^{4k+2})_\infty\label{krr10}.
\end{eqnarray}
\end{theo}
\begin{proof}
 For identities (\ref{krr1})-(\ref{krr4}), first set  $(a,b)=(0,0)$,
  $(-q^{-1},0)$, $(-q^{-1/2},0)$ and $(-q^{-1/2},-q^{-1})$ in
(\ref{zq2}), respectively,  and then apply (\ref{Ja1}).

For identities (\ref{krr9})-(\ref{krr11}), first set $(a,b)=(0,0)$, $(-q^{-1/2},0)$
and $(-q^{-1},0)$ in (\ref{zq}), respectively,  and then apply
(\ref{Ja2}).

Note that, for (\ref{krr3}), (\ref{krr4}) and
(\ref{krr10}),  we need to replace $q$ by $q^2$ at last.
\end{proof}

\begin{theo} For $k\geq 1$, the following identities hold
\begin{eqnarray}
&&\hskip -6 cm \sum_{l(\la)\leq
k}\frac{q^{|\la|+n_2(\la)-(\la_1^2+3\la_1)/2}(-q)_{\la_1}
(1-q^{\la_1})}{(-q)_{\la_k}(q)_\la}\nonumber\\
&=&\frac{(-q)_\infty}{(q)_\infty}(q^{2k+1},q^{2k-1},q^2;q^{2k+1})_\infty\label{krr5},\\
&&\hskip -6 cm \sum_{l(\la)\leq k}\frac{q^{2|\la|+2n_2(\la)-\la_1^2-2\la_1}
(-q;q^2)_{\la_1}}{(-q^2;q^2)_{\la_k}(q^2;q^2)_\la}\nonumber\\
&=&\frac{(-q;q^2)_\infty}{(q^2;q^2)_\infty}(q^{4k+2},q^{4k-1},q^3;q^{4k+2})_\infty\label{krr7},\\
\sum_{l(\la)\leq
k}\frac{q^{|\la|+n_2(\la)-\la_1}}{(-q)_{\la_k}(q)_\la}
&=&\frac{(q^{2k+2},q^{2k},q^2;q^{2k+2})_\infty}{(q)_\infty}\label{krr6},\\
\sum_{l(\la)\leq
k}\frac{q^{|\la|+n_2(\la)-2\la_1}(1-q^{2\la_1})}{(-q)_{\la_k}(q)_\la}
&=&\frac{(q^{2k+2},q^{2k-1},q^3;q^{2k+2})_\infty}{(q)_\infty}\label{krr8},\\
\sum_{l(\la)\leq k}\frac{q^{n_2(\la)-\la_1}}{(-q)_{\la_k}(q)_\la}
&=&\frac{(-1)_\infty}{(q^2;q^2)_\infty}(q^{2k+2},q^{k},q^{k+2};q^{2k+2})_\infty\label{krr12}.
\end{eqnarray}
\end{theo}
\begin{proof}
For $i\in\{0,1,2\}$, denote by $[b^i]$ the operation of extracting
the coefficient of $b^i$ in the corresponding identity. For
(\ref{krr5})-(\ref{krr8}), apply  the following operations to
(\ref{zq2}) respectively :  $a=-q^{-1}$ and $(1-1/q)[b]$, $a=0$
and $[b^0]+(1-1/q)[b]$, $a=-q^{-1/2}$ and $[b^0]+(1-1/q)[b]$,
$a=0$ and $[b]+(1-1/q)[b^2]$, and then apply (\ref{Ja1}).
Note that, for (\ref{krr7}), we need to replace  $q$ by $q^2$  at last.

 For (\ref{krr12})  apply the operations $a=0$
and $(1-1/q)[b]$ to (\ref{zq}) and then apply (\ref{Ja2}).
\end{proof}

\begin{rem}
As speculated by the anonymous referee, all of the Rogers-Ramanujan type identities given in Theorems~3 and 4 are known.

For example, specializing
equation (3.4) of
Bressoud~\cite{Br2}  (see also \cite{Br1}) with $k\rightarrow k+1$ and $r=1$, $r=2$ and $r=k+1$, respectively,
 we recover identities (\ref{krr1}), (\ref{krr6})) and  (\ref{krr9});
while specializing equation (3.9) of Bressoud~\cite{Br2}
with $k\rightarrow k+1$ and $r=1$, $r=2$ and $r=k+1$, respectively, we recover
  \eqref{krr3}, \eqref{krr7} and (\ref{krr10}).

Since we derived  all
these identities in Theorems~3 and 4 from
the two master identities (13) and (14),
instead of identifying each identity individually,
it suffices to
identify the later two with known results in the literature.
In 1984,
by means of Bailey chains, Andrews proved a
remarkable generalization of Bailey's lemma
\cite[Thm. 1]{An84}, which contains many multiple Rogers-Ramanujan type identities as special cases.
In particular, identities (13) and (14) are limit cases of Andrews' theorem.
More precisely,  to derive (13), set $a=q$, $b_k=1/a$, $c_k=1/b$ in
Andrews' formula and
let $N, b_1,c_1,\ldots, b_{k-1}, c_{k-1}\to \infty$, finally apply
the Bailey pair  E(3) of Slater's paper~\cite{Sl1}.
To derive (14) we do the same thing except that we set $a=1$ and apply
 the Bailey pair B(3) of Slater's paper~\cite{Sl1}.

\end{rem}

When $k=1$, identities (\ref{krr1}), (\ref{krr2}), (\ref{krr4}),
(\ref{krr8}) and (\ref{krr11}) reduce directly to special cases of the
$q$-binomial identity (\ref{qbinomial}). For example, when $k=1$
identity (\ref{krr4}) reduces to
$$
\sum_{n=0}^\infty\frac{q^n(-q;q^2)_n}{(q^2;q^2)_n}\frac{(-q^2;q^2)_\infty}{(q;q^2)_\infty},
$$
which is the $q$-binomial identity~(\ref{qbinomial}) after
substitutions $q\rightarrow q^2$, $a\rightarrow -q$ and $x\rightarrow
q$. The other identities reduce to the following
Rogers-Ramanujan type identities :
\begin{eqnarray}
\sum_{n=0}^\infty{q^{n^2+2n}(-q;q^2)_n\over (q^4;q^4)_{n}}
&=&\frac{(-q;q^2)_\infty}{(q^2;q^2)_\infty}(q,q^5,q^6;q^6)_\infty\label{rr3},\\
\sum_{n=0}^\infty{q^{n^2}\over (q^2;q^2)_{n}}&=&\frac{(q^2,q^2,q^4;q^4)_\infty}{(q)_\infty}\label{rr6},\\
\sum_{n=0}^\infty{q^{n^2}(-q;q^2)_n\over (q^4;q^4)_{n}}&=&\frac{(-q;q^2)_\infty}{(q^2;q^2)_\infty}(q^3,q^3,q^6;q^6)_\infty\label{rr7}.
\end{eqnarray}
Note that (\ref{rr6}) is again a special case of the $q$-binomial identity (1)
and (\ref{rr7}) is  (25)
of Slater's list~\cite{Sl2}.

\section{Proof of Theorem~1 and consequences}
\subsection{Proof of Theorem~1}
For any statement $A$ it will be convenient to introduce the Boolean function
 $\chi(A)$, which is 1 if $A$ is true
 and 0 if $A$ is false.
Consider the generating function $$
S(u)=\sum_{\la_0,\la}\left(\prod_{i=1}^{\la_0-1}(-q)_{m_i}\right)P_\la(X,q^2)\,
u^{\la_0}
$$ where the sum is over all partitions $\la=(\la_1,\ldots,
\la_n)$ and the integers
 $\la_0\geq \la_1$. Suppose
$\la=(\mu_1^{r_1}\, \mu_2^{r_2}\, \ldots \mu_k^{r_k})$,
 where
$\mu_1>\mu_2>\cdots >\mu_k\geq 0$ and $(r_1, \ldots, r_k)$ is a
 composition of $n$.

Let $S_n^\la$ be the set of permutations of $S_n$ which fix $\la$.
 Each $w\in S_n/S_n^\la$ corresponds to  a surjective mapping $f:
X\longrightarrow \{1,2,\ldots, k\}$ such that $|f^{-1}(i)|=r_i$.
For any subset $Y$ of $X$, let $p(Y)$ denote  the product of the
elements of $Y$ (in particular, $p(\emptyset)=1$). We can rewrite
Hall-Littlewood functions as follows : $$ P_\la(X,q^2)=\sum_{f}
p(f^{-1}(1))^{\mu_1}\cdots p(f^{-1}(k))^{\mu_k}
\prod_{f(x_i)<f(x_j)}{x_i-q^2x_j\over x_i-x_j}, $$ summed over all
surjective mappings $f: X\longrightarrow \{1,2,\ldots, k\}$ such
that $|f^{-1}(i)|=r_i$. Furthermore, each such $f$ determines a
\emph{filtration}  of $X$ :
\begin{equation}\label{filtration}
 {\cal F}: \quad
\emptyset=F_0\subsetneq F_1\subsetneq \cdots \subsetneq F_k=X,
\end{equation}
according to the rule $x_i\in F_l\Longleftrightarrow f(x_i)\leq l$
for $1\leq l\leq k$. Conversely, such a filtration ${\cal
F}=(F_0,\, F_1, \ldots, F_k)$ determines a surjection $f:
X\longrightarrow \{1,2,\ldots, k\}$ uniquely. Thus we can write :
\begin{equation}\label{filtre}
P_\la(X,q^2)=\sum_{\cal F}\pi_{\cal F}\prod_{1\leq i\leq
k}p(F_i\setminus F_{i-1})^{\mu_i},
\end{equation}
summed  over all the filtrations $\cal F$  such that
 $|F_i|=r_1+r_2+\cdots +r_i$ for $1\leq i\leq k$, and
$$ \pi_{\cal F}=\prod_{f(x_i)<f(x_j)}{x_i-q^2x_j\over x_i-x_j}, $$
where $f$ is the function defined by $\cal F$.

Now let $\nu_i=\mu_i-\mu_{i+1}$ if $1\leq i\leq k-1$ and
$\nu_k=\mu_k$, thus $\nu_i>0$ if $i<k$ and $\nu_k\geq 0$. Furthermore, let  $\mu_0=\la_0$ and
$\nu_0=\mu_0-\mu_1$ in the definition of $S(u)$, so that
$\nu_0\geq 0$ and $\mu_0=\nu_0+\nu_1+\cdots +\nu_k$. Define
 $c_{\cal F}=\prod_{i=1}^k(-q)_{|F_i\setminus
F_{i-1}|}$ for any filtration $\cal F$. Thus, since the  lengths of columns of $\la$ are
 $|F_j|=r_1+\cdots +r_j$
with multiplicities $\nu_j$  and $r_j=m_{\mu_j}(\la)$ for $1\leq j\leq k$, we have
\begin{eqnarray*}
\prod_{i=1}^{\la_0-1}(-q)_{m_i}&=&c_{\cal F}\times\left(\chi(\nu_k=0)
(-q)_{|F_k\setminus F_{k-1}|}+\chi(\nu_k\neq 0)\right)^{-1}\\
&&\hskip 0.7cm\times\left(\chi(\nu_0=0)(-q)_{|F_1|}+\chi(\nu_0\neq 0)\right)^{-1}.
\end{eqnarray*}
Let $F(X)$ be the set of filtrations of $X$. Summarizing we obtain
\begin{eqnarray}
S(u)&=&\sum_{{\cal F}\in F(X)}c_{\cal F}\,\pi_{\cal F}\,
\sum_{\nu_1>0}(u\,p(F_1))^{\nu_1}\cdots\sum_{\nu_{k-1}>0}(u\,p(F_{k-1}))^{\nu_{k-1}}\nonumber\\
&&\hskip 1 cm\times\sum_{\nu_0\geq 0}{u^{\nu_0}\over
\chi(\nu_0=0)\,(-q)_{|F_1|} +\chi(\nu_0\neq 0)}\nonumber\\
&&\hskip 1 cm\times\sum_{\nu_k\geq 0} {u^{\nu_k}\,p(F_k)^{\nu_k}\over
\chi(\nu_k=0)\,(-q)_{|F_k\setminus F_{k-1}|}+\chi(\nu_k\neq
0)}\label{bigmac'}.
\end{eqnarray}
For any filtration $\cal F$ of $X$ set
$$ {\cal A}_{\cal F}(X,
u)=c_{\cal F}\,\prod_{j}\left[{p(F_j)u \over 1-p(F_j)u}+{\chi(F_j=X)\over
(-q)_{|F_j\setminus F_{j-1}|}}+{\chi(F_j=\emptyset)\over
(-q)_{|F_1|}}\right].
$$
It follows from (\ref{bigmac'}) that
$$
S(u)=\sum_{{\cal F}\in F(X)}\pi_{\cal F}{\cal A}_{\cal F}(X, u).
$$
 Hence  $S(u)$ is a rational function of $u$ with simple
 poles at  $1/p(Y)$, where $Y$ is a subset of $ X$.
 We are
now proceeding to compute  the corresponding  residue $c(Y)$
 at each pole $u=1/p(Y)$.

Let us  start with $c(\emptyset)$. Writing $\la_0=\la_1+k$ with
$k\geq 0$, we see that
\begin{eqnarray*}
S(u)&=&\sum_\la f_{\la}(q)P_\la(X,q^2)u^{\la_1}\sum_{k\geq 0}{u^k\over
\chi(k=0)(-q)_{m_{\la_1}}+\chi(k\neq 0)}\\
 &=&\sum_\la f_{\la}(q)P_\la(X,q^2)u^{\la_1}\left({u\over 1-u}+{1\over
(-q)_{m_{\la_1}}}\right).
\end{eqnarray*}
It follows from (\ref{kawa1}) that
$$
c(\emptyset)=\left[S(u)(1-u)\right]_{u=1}=\Phi(X).
$$
For the computations of other  residues, we need some more
notations. For any $Y\subseteq X$, let $Y'=X\setminus Y$ and
$-Y=\{x_i^{-1}:x_i\in Y\}$. Then
\begin{equation}\label{residu}
c(Y)=\left[\sum_{\cal F}\pi_{\cal F}{\cal A}_{\cal
F}(X,\,u)(1-p(Y)u)\right]_{u=p(-Y)}.
\end{equation}
If $Y\notin\cal F$, the corresponding summand is  equal to  0.
Thus we need only to consider the following filtrations ${\cal F}$~:
 $$ \emptyset=F_0\subsetneq \cdots \subsetneq F_t=Y\subsetneq
\cdots \subsetneq  F_k=X\qquad 1\leq t\leq k.
$$
We may then split
$\cal F$ into two filtrations ${\cal F}_1$ and ${\cal F}_2$~ :
\begin{eqnarray*}
{\cal F}_1&:& \emptyset \subsetneq  -(Y\setminus
F_{t-1})\subsetneq \cdots \subsetneq -(Y\setminus F_1)\subsetneq
-Y,\\ {\cal F}_2&:& \emptyset \subsetneq F_{t+1}\setminus
Y\subsetneq   \cdots \subsetneq F_{k-1}\setminus Y\subsetneq  Y'.
\end{eqnarray*}
Then, writing  $v=p(Y)u$ and $c_{\cal F}=c_{{\cal F}_1}\times
c_{{\cal F}_2}$, we have
 $$
\pi_{\cal F}(X)=\pi_{{\cal F}_1}(-Y)\pi_{{\cal
F}_2}(Y')\prod_{x_i\in Y, x_j\in
Y'}\frac{1-q^2x_i^{-1}x_j}{1-x_i^{-1}x_j},
$$ and ${\cal A}_{\cal
F}(X,\,u)(1-p(Y)u)$ is equal to
\begin{eqnarray*}
&&{\cal A}_{{\cal F}_1}(-Y,\,v) {\cal A}_{{\cal
F}_2}(Y',\,v)(1-v)\left(\frac{v}{1-v}+\frac{\chi(Y=X)}
{(-q)_{|Y\setminus F_{t-1}|}}\right)\\
&&\times\left(\frac{v}{1-v}+
\frac{1}{(-q)_{|Y\setminus F_{t-1}|}}\right)^{-1}
\left(\frac{v}{1-v}
+\frac{1}{(-q)_{|F_{t+1}\setminus Y|}}\right)^{-1}.
\end{eqnarray*}
Thus when $u=p(-Y)$, i.e., $v=1$,
\begin{eqnarray*}
&&\hskip -1cm\left[\pi_{\cal F}(X){\cal A}_{\cal
F}(X,\,u)(1-p(Y)u)\right]_{u=p(-Y)}=\\ &&\left[\pi_{{\cal F}_1}(-Y){\cal
A}_{{\cal F}_1}(-Y,\,v)(1-v)\pi_{{\cal F}_2}(Y') {\cal A}_{{\cal
F}_2}(Y',\,v)(1-v)\right]_{v=1}\\ &&\hskip 2cm\times\prod_{x_i\in Y, x_j\in
Y'}\frac{1-q^2x_i^{-1}x_j}{1-x_i^{-1}x_j}.
\end{eqnarray*}
Using (\ref{residu}) and the result of $c(\emptyset)$, which can
be written
$$
\left[\sum_{\cal F}\pi_{\cal F}{\cal A}_{\cal F}(X,
u)(1-u)\right]_{u=1}=\Phi(X),
$$
 we get
$$
c(Y)=\Phi(-Y)\Phi(Y')\prod_{x_i\in Y, x_j\in
Y'}\frac{1-q^2x_i^{-1}x_j}{1-x_i^{-1}x_j}. $$ Each subset $Y$ of
$X$ can be encoded by  a sequence $\xi\in \{\pm 1\}^n$ according
to the rule~: $\xi_i=1$ if $x_i\notin Y$ and $\xi_i=-1$ if $x_i\in
Y$. Hence
$$ c(Y)=\Phi(X^{\xi}).
$$
Note also that
$$ p(Y)=\prod_i
x_i^{(1-\xi_i)/2},\qquad p(-Y)=\prod_i x_i^{(\xi_i-1)/2}.
$$
Now, extracting the coefficients of $u^k$ in the equation~:
$$ S(u)=\sum_{Y\subseteq X} {c(Y)\over 1-p(Y)u},$$
yields $$ \sum_{\la_1\leq
k}\left(\prod_{i=1}^{k-1}(-q)_{m_i}\right)P_{\la}(X,q^2)=\sum_{Y\subseteq
X}c(Y)p(Y)^k.
$$
Finally, substituting  the value of $c(Y)$ in the above formula
we obtain (\ref{extkawa1}).

\subsection{Proof of Corollary 1}
Recall \cite[p. 213]{Ma} that if $x_i=zq^{2i-2}$ ($1\leq
i\leq n$) then~:
\begin{equation}\label{value1}
P_{\la'}(X,q^2)=z^{|\la|}q^{2n(\la)}\left[{n\atop \la}\right]_{q^2}.
\end{equation}
Replacing each partition $\lambda$ by its  conjugate $\lambda'$ on the
left-hand side of (\ref{extkawa1})
yields the left-hand side of
(\ref{csqextkawa1}). Set
$$
\Psi(X)=\prod_{i}\frac{1}{1-x_i^2}\prod_{j<k}\frac{1-q^2x_jx_k}{1-x_jx_k}.
$$
Then, for any $\xi\in \{\pm 1\}^n$ such that the number of
$\xi_i=-1$ is $r$, $0\leq r\leq n$, we can write $\Phi(X^{\xi})$
as follows:
\begin{equation}\label{cal2}
\Phi(X^{\xi})=\Psi(X^{\xi})\,\prod_i\frac{1+qx_i^{\xi_i}}{1-x_i^{\xi_i}}(1-x_i^{2\xi_i}),
\end{equation}
which is readily seen to equal 0 unless $\xi\in \{-1\}^r\times
\{1\}^{n-r}$. Now, in the latter case,
 we have $\prod_ix_i^{k(1-\xi_i)/2}=z^{kr}q^{2k\left({r
\atop 2}\right)}$,
\begin{equation}\label{cal3}
\prod_{i=1}^n\frac{1+qx_i^{\xi_i}}{1-x_i^{\xi_i}}(1-x_i^{2\xi_i})=\frac{(z^2;q^4)_{n}}{z^{2r}q^{4\left({r \atop
2}\right)-r}}\frac{(-z/q;q^2)_r}{(z;q^2)_r}\frac{(-zq^{2r+1};q^2)_{n-r}}{(zq^{2r};q^2)_{n-r}},
\end{equation}
and \cite[p. 476]{St}~:
\begin{equation}\label{cal4}
\Psi(X^{\xi})=(-1)^rz^{2r}q^{6\left({r \atop
2}\right)}\left[{n\atop r}\right]{1-z^2q^{4r-2}\over
(zq^{r-1})_{n+1}}.
\end{equation}
Substituting these into the right-hand side of (\ref{extkawa1}) we obtain the
right-hand side of (\ref{csqextkawa1}) after simple manipulations.

When $n\to +\infty$, since
$\qbi{n}{\la}{}\to \frac{1}{(q)_\la}$,
equation (\ref{csqextkawa1}) reduces to~:
\begin{equation}\label{lim1}
\sum_{l(\la)\leq k}\frac{z^{|\la|}q^{2n(\la)}}
{(-q)_{\la_k}(q)_\la} =(-z/q)_\infty \sum_{r\geq
0}{(-1)^rz^{kr}q^{r+(2k+2)\bi{r}{2}}\over
(q^2;q^2)_{r}(z^2q^{2r-2})_\infty}(1-zq^{2r-1}).
\end{equation}
Furthermore, as in Section~2, setting
$z=q^2$ and $z=q$ in (\ref{lim1}) yields  (\ref{krr1}) and (\ref{krr9}), respectively.

\section{Elementary approach and proof of Theorem 2}
\subsection{Preliminaries}
We will need the following result, which corresponds to the case
$k\to\infty$ in  (\ref{csqextkawa1}), and can be proved in an
elementary way :
\begin{lem} For $n\geq 0$
\begin{equation}\label{kinf}
\sum_{\la}z^{|\la|}q^{2n(\la)}(-q)_\la\left[{n\atop \la}\right]_{q^2}=\frac{(-z)_{2n}}{(z^2;q^2)_n}.
\end{equation}
\end{lem}
\begin{proof}
Recall the following identity, which is proved in \cite{JZ2} :
\begin{equation}\label{qpieri}
q^{\bi{m}{2}+n(\mu)}\qbi{n}{m}{}\qbi{n}{\mu}{}=\sum_{\la}q^{n(\la)}\qbi{n}{\la}{}
\prod_{i\geq 1}\qbi{\la_i-\la_{i+1}}{\la_i-\mu_i}{},
\end{equation}
where the sum is over all partitions $\la$ such that $\la/\mu$ is
a horizontal $m$-strip, i.e., $\mu\subseteq\la$, $|\la/\mu|=m$
and there is at most one cell in each column of the Ferrers
diagram of $\la/\mu$.\\
We also need
\begin{equation}\label{mac2}
\sum_\la
z^{|\la|}q^{n(\la)}\qbi{n}{\la}{}=\frac{(-z)_n}{(z^2)_n},
\end{equation}
which can be found in \cite{JZ2, St}.\\
Using (\ref{mac2}) with $q$ replaced  by $q^2$ and (\ref{qbin}), the right-hand side of (\ref{kinf}) can be written
\begin{eqnarray*}
\frac{(-z;q^2)_n}{(z^2;q^2)_n}(-zq;q^2)_n&=&\sum_{\mu,\,m}
z^{|\mu|}q^{2n(\mu)}\qbi{n}{\mu}{q^2}z^{m}q^{2\left({m\atop
2}\right)+m}\qbi{n}{m}{q^2}\\
&=&\sum_{\la,\,m}z^{|\la|}q^{2n(\la)}\qbi{n}{\la}{q^2} \prod_{i\geq
1}\sum_{r_i\geq 0}q^{r_i}\qbi{\la_i-\la_{i+1}}{r_i}{q^2},
\end{eqnarray*}
where the last equality follows from (\ref{qpieri}), setting $r_i=\la_i-\mu_i$ for $i\geq 1$. Now we conclude by using (\ref{csqeuler}).
\end{proof}

Recall the following extension of the $n\to\infty$ case of (\ref{kinf}), which is Stembridge's lemma 3.3 (b) in \cite{St}, and identity (60) in \cite{JZ2} :
\begin{equation}\label{abext}
\sum_{\la}z^{|\la|}q^{2n(\la)}\frac{(a,b;q^{-1})_{\la_1}}{(q)_\la}=\frac{(az,bz)_{\infty}}{(z,abz)_\infty}.
\end{equation}
Now, using (\ref{kinf}), we are able to
prove directly identity (\ref{csqextkawa1}) in Corollary~1,
and then using (\ref{abext}),
to deduce an elementary proof of (\ref{ab}) in Theorem~2.

\subsection{Elementary proof of Corollary~1}
  Consider the
generating function of the left-hand side of (\ref{csqextkawa1})~:
\begin{eqnarray}
\varphi(u) &=&\sum_{k\geq 0}u^k\sum_{l(\la)\leq
k}\frac{(-q)_{\la}}{(-q)_{\la_k}}z^{|\la|}q^{2n(\la)}\left[{n\atop \la}\right]_{q^2}\\ &=&\sum_{\la}
u^{l(\la)}z^{|\la|}q^{2n(\la)}(-q)_\la \left[{n\atop
\la}\right]_{q^2}\sum_{k\geq 0} {u^k\over
(-q)_{\la_{k+l(\la)}}}\nonumber\\
&=&\sum_{\la}u^{l(\la)}z^{|\la|}q^{2n(\la)}(-q)_\la
\left[{n\atop
\la}\right]_{q^2}\left(\frac{u}{1-u}+\frac{1}{(-q)_{\la_{l(\la)}}}\right)
\label{ref1},
\end{eqnarray}
where the last equality follows from the fact that $\displaystyle\la_{k+l(\la)}=0$ unless $k=0$.
Now, each partition $\la$ with parts bounded by $n$ can be encoded
by
 a pair of sequences
 $\nu=(\nu_0, \nu_1, \cdots, \nu_l)$ and ${\mathbf m}=(m_0,\cdots,m_{l})$
such that $\la=(\nu_0^{m_0},\ldots, \nu_l^{m_l})$, where
$n=\nu_0>\nu_1> \cdots
>\nu_l>0$ and $\nu_i$ has multiplicity $m_i\geq 1$ for $1\leq
i\leq l$ and $\nu_0=n$ has   multiplicity $m_0\geq 0$. Using the
notation~: $$ <\alpha>=\frac{\alpha}{1-\alpha},\quad
u_i=z^iq^{i(i-1)}\quad\hbox{for}\quad
 i\geq 0,
$$ we can then rewrite (\ref{ref1}) as follows~:
\begin{eqnarray}
\varphi(u)&=&\sum_{\nu}(-q)_\nu\left[{n\atop
\nu}\right]_{q^2}\left(<u>+\frac{1}{(-q)_{\nu_l }}\right)\nonumber\\
&&\hskip 10pt \times\sum_{{\mathbf
m}}\left((u_nu)^{m_0}+{\chi(m_0=0)\over
(-q)_{n-\nu_1}}\right)\prod_{i=1}^{l}(u_{\nu_i}u)^{m_i
}\nonumber\\ &=&\sum_{\nu}{(q^2;q^2)_{n}\over
(q)_\nu}B_\nu,\label{x}
\end{eqnarray}
where the sum is over all strict partitions
$\nu=(\nu_0,\nu_1,\ldots, \nu_l)$ and $$
B_\nu=\left(<u>+\frac{1}{(-q)_{\nu_l }}\right)\left(
<u_ru>+\frac{1}{(-q)_{n-\nu_1}}\right)
\prod_{i=1}^{l}<u_{\nu_i}u>. $$ So $\varphi(u)$ is  a rational
fraction with simple poles at $u_r^{-1}$ for $0\leq r\leq n$. Let
$b_r(z,n)$ be the corresponding residue of
 $\varphi(u)$ at $u_r^{-1}$ for $0\leq r\leq n$.
Then, it follows from (\ref{x}) that
\begin{equation}\label{residue}
b_r(z,n)=\sum_\nu{(q^2;q^2)_{n}\over
(q)_\nu}\left[B_\nu(1-u_ru)\right]_{u=u_r^{-1}}.
\end{equation}
We shall first consider the cases where $r=0$ or $n$. Using
 (\ref{ref1}) and (\ref{kinf}) we have
\begin{equation}\label{b0}
b_0(z,n)=\left[\varphi(u)(1-u)\right]_{u=1}=\frac{(-z)_{2n}}{(z^2;q^2)_n}.
\end{equation}
Now, by (\ref{x}) and(\ref{residue}) we have
\begin{equation}\label{b0bis}
b_0(z,n)=\sum_{\nu}{(q^2;q^2)_{n}\over (q)_\nu}
\left(<u_n>+\frac{1}{(-q)_{n-\nu_1}}\right)
\prod_{i=1}^{l}<u_{\nu_i}>,
\end{equation}
and
\begin{equation}\label{refb}
b_n(z,n)=\sum_{\nu}{(q^2;q^2)_{n}\over
(q)_\nu}\left(<1/u_n>+\frac{1}{(-q)_{\nu_l }}\right)
\prod_{i=1}^{l}<u_{\nu_i}/u_n>,
\end{equation}
 which, by setting
$\mu_i=n-\nu_{l+1-i}$ for $1\leq i\leq l$ and $\mu_0=n$,   can be
written as
\begin{equation}\label{br}
b_n(z,n)=\sum_{\mu}{(q^2;q^2)_{n}\over
(q)_\mu}\left(<1/u_n>+\frac{1}{(-q)_{n-\mu_1
}}\right)\prod_{i=1}^{l}<u_{n-\mu_i}/u_n>.
\end{equation}
Comparing (\ref{br}) with (\ref{b0bis}) we see that $b_n(z,n)$ is
equal to $b_0(z,n)$ with $z$ replaced by $z^{-1}q^{-2n+2}$. It
follows from (\ref{b0}) that
\begin{equation}\label{ref3}
b_n(z,n)=b_0(z^{-1}q^{-2n+2},n)=(-1)^nq^{n^2}\frac{(-z/q)_{2n}}{(z^2q^{2n-2};q^2)_n}.
\end{equation}
Consider now the case where  $0<r<n$. Clearly, for each partition
$\nu$, the corresponding
 summand in (\ref{residue}) is not zero only if
$\nu_j=r$ for some $j$, $0\leq j\leq n$. Furthermore, each such
partition $\nu$ can be split into two strict partitions
$\rho=(\rho_0, \rho_1,\ldots, \rho_{j-1})$ and $\sigma=(\sigma_0,
 \ldots, \sigma_{l-j})$ such that
$\rho_i=\nu_i-r$ for $0\leq i\leq j-1$ and $ \sigma_s=\nu_{j+s}$
for $0\leq s\leq l-j$. So we can write (\ref{residue}) as
follows~:
\begin{eqnarray*}
b_r(z, n)&=&\left[{n\atop r}\right]_{q^2}
\sum_{\rho}\frac{(q^2;q^2)_{n-r}}{(q)_\rho}F_\rho(r)\times\sum_\sigma\frac{(q^2;q^2)_{r}}{(q)_\sigma}G_\sigma(r)
\end{eqnarray*}
where  for $\rho=(\rho_0, \rho_1,\ldots, \rho_l)$ with
$\rho_0=n-r$, $$ F_\rho(r)=\left(<u_n/u_r>+{1\over
(-q)_{n-r-\rho_1}}\right)\prod_{i=1}^{l(\rho)}<u_{\rho_i+r}/u_r>,
$$ and for $\sigma=(\sigma_0, \ldots, \sigma_l)$ with
$\sigma_0=n$, $$ G_\sigma(r)=\left(<1/u_r>
+\frac{1}{(-q)_{\sigma_{l} }}\right)
\prod_{i=1}^{l(\sigma)}<u_{\sigma_i}/u_r>. $$ Comparing with
(\ref{b0bis}) and (\ref{br}) and using (\ref{b0}) and (\ref{ref3})
we obtain
\begin{eqnarray*}
b_r(z,n)&=&\left[{n\atop r}\right]_{q^2}b_0(zq^{2r},n-r)\,b_r(z,r)\\
 &=&(-1)^rq^{r+2\left({n\atop
r}\right)}\frac{(-z/q)_{2n+1}}{(z^2q^{2r-2},q^2)_{n+1}}(1-zq^{4r-1}).
\end{eqnarray*}
Finally,  extracting the coefficients of $u^k$ in the equation $$
\varphi(u)=\sum_{p=0}^n\frac{b_r(z,n)}{1-u_ru}, $$ and using  the
values for $b_r(z,n)$ we obtain (\ref{csqextkawa1}).

\subsection{Proof of Theorem~2}
Consider the generating function of the left-hand side of
(\ref{ab}) :
\begin{eqnarray}
\varphi_{ab}(u)&:=&\sum_{k\geq 0}u^k\sum_{l(\la)\leq
k}z^{|\la|}q^{2n(\la)}\frac{(a,\,
b;q^{-1})_{\la_1}}{(q)_\la(-q)_{\la_k}}\nonumber\\
&=&\sum_\la\sum_{k\geq
0}u^{k+l(\la)}z^{|\la|}q^{2n(\la)}\frac{(a,\,
b;q^{-1})_{\la_1}}{(q)_\la(-q)_{\la_{l(\la)+k}}},\nonumber
\end{eqnarray}
where the sum is over all the partitions $\lambda$, and as before $\la_{l(\la)+k}=0$ unless $k=0$. Thus
\begin{equation}\label{g}
\varphi_{ab}(u)
=\sum_\la u^{l(\la)}z^{|\la|}q^{2n(\la)}\frac{(a,\,
b;q^{-1})_{\la_1}}{(q)_\la}\left(\frac{u}{1-u}+\frac{1}{(-q)_{\la_{
l(\la)}}}\right).
\end{equation}
As in the elementary proof of Corollary~1, we can replace
any  partition $\la$  by a pair $(\nu, {\bf m})$, where
$\nu$ is a strict partition consisting of distinct parts
   $\nu_1, \cdots, \nu_l$ of $\la$, so that $\nu_1> \cdots
>\nu_l>0$,  and ${\bf m}=(m_1,\ldots, m_l)$
is the sequence of multiplicities of
$\nu_i$ for $1\leq i\leq l$. Therefore
\begin{eqnarray}
\varphi_{ab}(u)&=&\sum_{\nu,\,{\mathbf m}}\frac{(a,\,
b;q^{-1})_{\nu_1}}{(q)_\nu}\left(\frac{u}{1-u}+\frac{1}{(-q)_{\nu_l
}}\right) \prod_{i=1}^{l}(u_{\nu_i}u)^{m_i}\nonumber\\
&=&\sum_{\nu}\frac{(a,\,
b;q^{-1})_{\nu_1}}{(q)_\nu}\left(<u>+\frac{1}{(-q)_{\nu_l}}\right)
\prod_{i=1}^l<u_{\nu_i}u>\label{h},
\end{eqnarray}
where the sum is over all the strict partitions $\nu$.
Each of the terms in this sum, as a rational function of $u$, has
a finite set of simple poles, which may occur at the points
$u_r^{-1}$ for $r\geq 0$. Therefore, each term is a linear
combination of partial fractions. Moreover, the sum of their
expansions converges coefficientwise. So $\varphi_{ab}$ has an
expansion
$$
\varphi_{ab}(u)=\sum_{r\geq
0}\frac{c_r}{1-uz^rq^{r(r-1)}},
$$
where $c_r$ denotes the formal
sum of partial fraction coefficients contributed by the terms of
(\ref{h}).
It remains to compute these residues $c_r$ $(r\geq 0)$.
By using (\ref{abext}) and (\ref{h}), we get immediately
$$
c_0=\left[\varphi_{ab}(u)(1-u)\right]_{u=1}= \frac{(az,\, bz)_\infty}{(z,\, abz)_\infty}.$$ In view of (\ref{h}),
this yields the identity
\begin{equation}\label{i}
\sum_{\nu}\frac{(a,\,
b;q^{-1})_{\nu_1}}{(q)_\nu}\prod_{i=1}^{l}
<u_{\nu_i}>=\frac{(az,\, bz)_\infty}{(z,abz)_\infty}.
\end{equation}
Clearly, a summand in (\ref{h}) has a non zero contribution to $c_r$ ($r>0$)
 only if the corresponding partition $\nu$ has a part equal to $r$.
 For any partition  $\nu$ such that  $\exists
j\, |\,\nu_j=r$,  set
$\rho_i:=\nu_i-r$ for $1\leq i<j$ and $\sigma_i:=\nu_{i+j}$ for
$0\leq i\leq l-j$, we then get two partitions $\rho$ and $\sigma$,  with
$\sigma_i$ bounded by $r$.
Multiplying (\ref{h}) by $(1-u_ru)$ and setting $u=1/u_r$  we obtain
\begin{eqnarray*}
c_r&=&\sum_{\rho}\frac{(a,b;q^{-1})_{\rho_1+r}}{(q)_\rho}
\prod_{i=1}^{j-1}<u_{r+\rho_i}/u_r>\\
&&\times\sum_{\sigma}\frac{1}{(q)_\sigma}
\left(<1/u_r>+\frac{1}{(-q)_{\sigma_{l-j}}}\right)
\prod_{i=1}^{l-j}<u_{\sigma_i}/u_r>.
\end{eqnarray*}
In view of (\ref{refb}) the inner sum over  $\sigma$
is equal to $b_r(z,r)/(q^2,q^2)_{r}$, and applying (\ref{ref3}) we get
\begin{eqnarray*}
c_r&=&  (-1)^rq^{r+2\left({r\atop
2}\right)}\frac{(-z/q)_{2r}}{(z^2q^{2r-2},q^2)_{r}}\frac{(a,b;q^{-1})_r}{(q^2;q^2)_{r}}\\
&&\hskip 2cm  \times\sum_{\rho}
\frac{(aq^{-r},bq^{-r};q^{-1})_{\rho_1}}
{(q)_\rho}\prod_{i=1}^{j-1}<{u_{r+\rho _i}/u_r}>.
\end{eqnarray*}
Now, the sum over $\rho$ can be computed using
 (\ref{i}) with $a$, $b$ and $z$
replaced  by $aq^{-r}$, $bq^{-r}$ and  $zq^{2r}$ respectively. After
simplification, we obtain
$$ c_r=(-1)^rq^{r+2\left({r\atop
2}\right)}\frac{(-z/q)_{\infty}}{(z^2q^{2r-2},q^2)_{\infty}}\frac{(a,b;q^{-1})_r}{(q^2;q^2)_{r}}\frac{(azq^r,\, bzq^r)_\infty}{(abz)_\infty},
$$
which completes the proof.

\section*{Acknowledgement} The second and third authors are
partially supported by EC's IHRP Programme, within the Research
Training Network ``Algebraic Combinatorics in Europe'', grant
HPRN-CT-2001-00272.



\small
\begin{thebibliography}{99}

\bibitem{An76} \textsc{Andrews} (G.E.), \emph{The theory of partitions},
Encyclopedia of mathematics and its applications,
 Vol. {\bf 2}, Addison-Wesley, Reading, Massachusetts, 1976.
\bibitem{An84} \textsc{Andrews} (G.E.),
Multiple series Rogers-Ramanujan type identities, \emph{Pacific J.
Math.}, Vol. {\bf 114}, No. 2, 267-283, 1984.
\bibitem{Br1} \textsc{Bressoud} (D. M.), An analytic generalization of the Rogers-Ramanujan identities with interpretation, \emph{Quart. J. Math. Oxford Ser. (2)}, {\bf 31}, no 124, 385-399, 1980.
\bibitem{Br2} \textsc{Bressoud} (D. M.), Analytic and combinatorial generalizations of the Rogers-Ramanujan identities, \emph{Mem. Amer. Math. Soc.}, {\bf 24},
no 227, 1980.

\bibitem{GR} \textsc{Gasper} (G.) and \textsc{Rahman} (M.),
\emph{ Basic Hypergeometric Series}, Second Edition, Encyclopedia of MathematicsAnd Its Applications 96, Cambridge University Press, Cambridge, 2004.


\bibitem{JZ1} \textsc{Jouhet} (F.) and \textsc{Zeng} (J.), Some
 new identities for Schur functions, \emph{Adv. Appl. Math.}, {\bf 27},
493-509, 2001.
\bibitem{JZ2} \textsc{Jouhet} (F.) and \textsc{Zeng} (J.), New
Identities of Hall-Littlewood Polynomials and  Applications, {\bf 10}, 89-112, 2005.

\bibitem{Kaw} \textsc{Kawanaka} (N.), A q-series identity involving Schur functions and related topics, \emph{Osaka J. Math.}, {\bf 36},
157-176, 1999.
\bibitem{Ma} \textsc{Macdonald} (I.G.),
\emph{Symmetric functions and Hall polynomials}, Clarendon Press,
second edition, Oxford, 1995.

\bibitem{Sl1} \textsc{Slater} (L. J.), A new  proof of Rogers's transformations of infinite series,
 \emph{Proc. London  Math. Soc.}, {\bf
53} (2), 460-475 (1951).
\bibitem{Sl2} \textsc{Slater} (L. J.), Further identities of
the  Rogers-Ramanujan Type, \emph{Proc. London  Math. Soc.}, {\bf
54} (2), 147-167 (1951-52).
\bibitem{St} \textsc{Stembridge} (J. R.), Hall-Littlewood
functions, plane partitions, and the Rogers-Ramanujan identities,
\emph{Trans. Amer. Math. Soc.}, {\bf 319}, no.2, 469-498, 1990.

\bibitem{Ze}
\textsc{Zeng} (J.), On the $q$-variations of Sylvester's bijection, \emph{The Ramanujan J.}, {\bf 9}, 289-303, 2005.


\end{thebibliography}
\end{document}